# Exact and simple formulas for the linearization coefficients of products of orthogonal polynomials and physical application

A. D. Alhaidari

*Saudi Center for Theoretical Physics, P.O. Box 32741, Jeddah 21438, Saudi Arabia*

**Abstract:** We obtain exact, simple and very compact expressions for the linearization coefficients of the products of orthogonal polynomials; both the conventional Clebsch-Gordan-type and the modified version. The expressions are general depending only on the coefficients of the three-term recursion relation of the linearizing polynomials. These are more appropriate and useful for doing numerical calculations when compared to other exact formulas found in the mathematics literature, some of which apply only to special class of polynomials while others may involve the evaluation of intractable integrals. As an application in physics, we present a remarkable phenomenon where nonlinear coupling in a physical system with pure continuous spectrum generates a mixed spectrum of continuous and discrete energies.

**Keywords**: linearization coefficients, orthogonal polynomials, Clebsch-Gordan-type coefficients, three-term recursion relation, discrete and continuous spectra.

## 1. Introduction

Orthogonal polynomials are exceptional elements in the set of special functions with a wide range of applications in various branches of science, engineering and economics [1-3]. The literature is populated with everlasting methods that utilize these objects in various scenarios and applications (see, for example, Refs. [4-7]). In physics, for example, orthogonal polynomials appear abundantly in quantum mechanics (mostly weighted by positive functions) to represent states of many interesting physical systems. In nonlinear physics, they appear as products due nonlinear couplings. Consequently, the study of such products of orthogonal polynomials and how to linearize theme become highly significant.

Given two sequences of orthogonal polynomials $\{p_n(x), q_n(x)\}$, the linearization expansion of the product $p_n(x) p_m(x)$ in terms of $q_n(x)$ reads as follows

$$p_n(x) p_m(x) = \sum_{k=0}^{n+m} c_{n,m}^k q_k(x), \qquad (1)$$

where $\{c_{n,m}^k\}$ are called the linearization coefficients. We consider two cases $q_n(x) = p_n(x)$ where $\{c_{n,m}^k\}$ are known as the Clebsch-Gordan-type (C-G-type) coefficients and $q_n(x) \neq p_n(x)$ where they are referred to as the modified C-G-type coefficients. This is an old problem in orthogonal polynomials whose solution is relevant across many fields in pure and applied mathematics and also in mathematical physics (see section 9 in Ref. [3]). Most formulas found in the literature for these coefficients, especially for the case $q_n(x) \neq p_n(x)$, are either very complicated, apply to special class of polynomials, or involve the evaluation of intractable integrals (see, for example Refs. [8-14]). Most, though, are not suitable for doing numerical calculations. The physically most interesting and highly non-trivial scenario is when the polynomials $\{p_n(x)\}$ have a pure continuous spectrum whereas the polynomials $\{q_n(x)\}$ are endowed with a mix of discrete and continuous spectrum. Physically, this implies that two physical systems with continuous energy spectra (no bound states) could nonlinearly couple to produce a mixed spectrum of continuous and discrete energies (i.e., generating bound states).

In Section 2, we consider the case $q_n(x) = p_n(x)$ whereas in Section 3, we study the case $q_n(x) \neq p_n(x)$. However, in both sections all polynomials have pure continuous spectra. In Section 4, we repeat the treatment of Section 3 but with $p_n(x)$ having a pure continuous spectrum while $q_n(x)$ having a mix of continuous and discrete spectrum. Finally in Section 5, we present a physical application where these findings are used to derive a solution for an interesting nonlinear scattering problem.

Throughout this study, we choose to work with the *orthonormal* version of polynomials where the orthogonality relation is normalized (i.e., $\langle p_n(x) | p_m(x) \rangle = \delta_{n,m}$) and the three-term recursion relation is symmetric. Nonetheless, for other versions (e.g., monic polynomials) our results can readily be used since all other versions are obtained as $p'_n(x) = T_n^p p_n(x)$, where $\{T_n^p\}$ are non-zero coefficients, giving $c'^k_{n,m} = (T_n^p T_m^p / T_k^q) c^k_{n,m}$.

## 2. Same polynomial species (C-G-type coefficients)

Let $\{p_n(x)\}$ be a set of orthonormal polynomials with pure continuous spectrum that satisfy the following symmetric three-term recursion relation

$$x p_n(x) = a_n p_n(x) + b_{n-1} p_{n-1}(x) + b_n p_{n+1}(x), \qquad (2)$$

for $n = 1, 2, \ldots$ and $b_n \neq 0$ for all $n$. The initial polynomials are $p_0(x) = 1$ and $p_1(x) = \gamma x + \eta$ such that $\gamma \neq 0$. They also satisfy the normalized orthogonality

$$\int_{x_-}^{x_+} \rho(x) p_n(x) p_m(x) dx = \delta_{n,m}, \qquad (3)$$

where $x_\pm$ are the boundaries of the space. These polynomials and their weight function $\rho(x)$ depend on the parameters $\gamma$, $\eta$ and any other parameters in the recursion coefficients $\{a_n, b_n\}$. In this setting, we can state the following Theorem:

*Theorem* **1**: In the linearization expansion

$$p_n(x) p_m(x) = \sum_{k=|n-m|}^{n+m} c^k_{n,m} p_k(x), \qquad (4)$$

the C-G-type coefficients are $c^k_{n,m} = [p_k(J)]_{n,m}$ where $J$ is the following tridiagonal symmetric matrix

$$J = \begin{pmatrix} a_0 & b_0 & & & \\ b_0 & a_1 & b_1 & & \\ & b_1 & a_2 & b_2 & \\ & & \times & \times & \times \\ & & & \times & \times \end{pmatrix}. \qquad (5)$$

***Proof***: Multiply both sides of Eq. (4) by $\rho(x)p_l(x)$ and integrate. Using the orthogonality (4) on the right side, one gets

$$\int_{x_-}^{x_+} \rho(x)p_l(x)p_n(x)p_m(x)dx = \sum_{k=|n-m|}^{n+m} c_{n,m}^k \delta_{l,k} = c_{n,m}^l. \tag{6}$$

Now, the recursion relation (2) and orthogonality (4) give

$$\int_{x_-}^{x_+} x\rho(x)p_n(x)p_m(x)dx = a_n\delta_{n,m} + b_{n-1}\delta_{n-1,m} + b_n\delta_{n+1,m} = J_{n,m}. \tag{7}$$

By iteration (i.e., replacing $x$ inside the integral by $x^2, x^3, ..., x^k$), we obtain

$$\int_{x_-}^{x_+} x^k \rho(x)p_n(x)p_m(x)dx = \left(J^k\right)_{n,m}. \tag{8}$$

Therefore, the left side of Eq. (6) becomes

$$\int_{x_-}^{x_+} p_l(x)\rho(x)p_n(x)p_m(x)dx = \left[p_l(J)\right]_{n,m}, \tag{9}$$

proving the theorem.  □

The expression for the C-G-type coefficients, $c_{n,m}^k = \left[p_k(J)\right]_{n,m}$, is compact and very simple, especially for doing numerical calculations. Moreover, it is a general expression that depends only on the recursion coefficients $\{a_n, b_n\}$. However, we can also obtain a recursion relation for $\{c_{n,m}^k\}$ so that we can evaluate all of them starting with some initial values. By multiplying both sides of Eq. (4) by $x$ and using the three-term recursion relation (2), we get

$$(a_n - a_k)c_{n,m}^k = b_k c_{n,m}^{k+1} - b_n c_{n+1,m}^k + b_{k-1}c_{n,m}^{k-1} - b_{n-1}c_{n-1,m}^k, \tag{10}$$

together with another recursion obtained from this one by the exchange $n \leftrightarrow m$ since $c_{n,m}^k = c_{m,n}^k$. Subtracting the two equations, we get the following four-term recursion relation for the same superscript $k$

$$(a_n - a_m)c_{n,m}^k = b_m c_{n,m+1}^k - b_n c_{m,n+1}^k + b_{m-1}c_{n,m-1}^k - b_{n-1}c_{m,n-1}^k, \tag{11a}$$

which can be solved for all $\{c_{n,m}^k\}_{n\neq m}$ starting with $\{c_{n,n}^k\}$ and $\{c_{n,n+1}^k\}$. This relation needs to be augmented by relation (10) with $n = m$ that reads

$$(a_n - a_k)c_{n,n}^k = b_k c_{n,n}^{k+1} + b_{k-1}c_{n,n}^{k-1} - b_n c_{n,n+1}^k - b_{n-1}c_{n,n-1}^k, \tag{11b}$$

which can be solved for all $\{c_{n,n+1}^k\}$ starting with $\{c_{n,n}^k\}$. Therefore, to solve for all $\{c_{n,m}^k\}$ using both (11a) and (11b), we need the seed values $\{c_{n,n}^k\}$.

Iterating (4) gives the linearization coefficients for the general product $p_n(x)p_m(x)p_k(x)...$
$...p_l(x)$. As an example, $p_n(x)p_m(x)p_k(x) = \sum_{l=0}^{n+m+k} C_{n,m,k}^l p_l(x)$, where $C_{n,m,k}^l = \sum_{j=|n-m|}^{n+m} c_{n,m}^j c_{j,k}^l$.

In the Appendix, we give two illustrative examples. The first is where $p_n(x)$ stands for the normalized Hermite polynomial and compare our findings to the well-known analytic results in the literature. In the second example, we take $p_n(x)$ as the normalized continuous dual Hahn polynomial with a pure continuous spectrum.

## 3. Dissimilar polynomial species (modified C-G-type coefficients)

In addition to the polynomials $\{p_n(x)\}$ above, let $\{q_n(x)\}$ be a different set of orthogonal polynomials but also with pure continuous spectrum and satisfying the following recursion relation

$$x q_n(x) = \alpha_n q_n(x) + \beta_{n-1} q_{n-1}(x) + \beta_n q_{n+1}(x), \tag{12}$$

for $n = 1, 2, ...$ and $\beta_n \neq 0$. The initial polynomials are $q_0(x) = 1$ and $q_1(x) = \kappa x + \zeta$ such that $\kappa \neq 0$. They satisfy the orthogonality

$$\int_{x_-}^{x_+} \omega(x) q_n(x) q_m(x) dx = \delta_{n,m}. \tag{13}$$

***Theorem 2***: In the linearization expansion

$$p_n(x) p_m(x) = \sum_{k=0}^{n+m} d_{n,m}^k q_k(x), \tag{14}$$

the modified C-G-type coefficients could be written either as $d_{n,m}^l = \sum_{k=|n-m|}^{n+m} c_{n,m}^k [p_k(K)]_{0,l}$

or as $d_{n,m}^l = [p_n(K) p_m(K)]_{0,l}$ where $K$ is identical to the matrix $J$ in (5) but with $\{a_n, b_n\} \mapsto \{\alpha_n, \beta_n\}$.

***Proof***: Using (4), we can rewrite (14) as

$$\sum_{k=|n-m|}^{n+m} c_{n,m}^k p_k(x) = \sum_{k=0}^{n+m} d_{n,m}^k q_k(x). \tag{15}$$

Using the orthogonality (13), we can write

$$d_{n,m}^l = \sum_{k=|n-m|}^{n+m} c_{n,m}^k \int_{x_-}^{x_+} \omega(x) p_k(x) q_l(x) dx. \tag{16}$$

To evaluate this integral, we follow the same procedure as in the previous section. Since $q_0(x) = 1$, we can rewrite the integral as follows

$$\int_{x_-}^{x_+} \omega(x) p_k(x) q_0(x) q_l(x) dx = \left[ p_k(K) \right]_{0,l}, \qquad (17)$$

Finally, we obtain

$$d_{n,m}^l = \sum_{k=|n-m|}^{n+m} c_{n,m}^k \left[ p_k(K) \right]_{0,l}. \qquad (18)$$

Like $c_{n,m}^k$, it is a general expression that depends only on the recursion coefficients $\{\alpha_n, \beta_n\}$ and $\{a_n, b_n\}$. An alternative, more compact and simpler but equivalent formula is obtained by multiplying the right sides of Eq. (14) by $\omega(x) q_l(x)$ and left side by $\omega(x) q_0(x) q_l(x)$ [since $q_0(x) = 1$] then integrating the equation to obtain

$$d_{n,m}^l = \left[ p_n(K) p_m(K) \right]_{0,l}. \qquad (19)$$

This expression is more compact and simpler. It depends only on the recursion coefficients $\{\alpha_n, \beta_n\}$. □

Similar to $\{c_{n,m}^k\}$, we can obtain a recursion relation for $\{d_{n,m}^k\}$ by multiplying both sides of Eq. (14) by $x$ then using the recursion relations (2) and (12) to obtain

$$(a_n - \alpha_k) d_{n,m}^k = \beta_{k-1} d_{n,m}^{k-1} - b_{n-1} d_{n-1,m}^k + \beta_k d_{n,m}^{k+1} - b_n d_{n+1,m}^k, \qquad (20)$$

together with another recursion obtained by $n \leftrightarrow m$ since $d_{n,m}^k = d_{m,n}^k$. Subtracting the two, we get the following four-term recursion relation for the same superscript $k$

$$(a_n - a_m) d_{n,m}^k = b_m d_{n,m+1}^k - b_n d_{m,n+1}^k + b_{m-1} d_{n,m-1}^k - b_{n-1} d_{m,n-1}^k, \qquad (21a)$$

which is identical to the recursion relation (10) for $c_{n,m}^k$. This relation needs to be supplemented by relation (20) with $n = m$ that reads

$$(a_n - \alpha_k) d_{n,n}^k = \beta_k d_{n,n}^{k+1} + \beta_{k-1} d_{n,n}^{k-1} - b_n d_{n,n+1}^k - b_{n-1} d_{n,n-1}^k. \qquad (21b)$$

Therefore, (21a) and (21b) can be solve for all $\{d_{n,m}^k\}$ starting with $\{d_{n,n}^k\}$.

Iterating (14) gives the linearization coefficients for the general product $p_n(x) p_m(x) p_k(x)...$ $...p_l(x)$. As an example, $p_n(x) p_m(x) p_k(x) = \sum_{l=0}^{n+m+k} D_{n,m,k}^l q_l(x)$, where $D_{n,m,k}^l = \sum_{j=|n-m|}^{n+m} c_{n,m}^j d_{j,k}^l$.

In the Appendix, we give two illustrative examples. The first is where $p_n(x)$ stands for the normalized Laguerre polynomial $L_n^\nu(x)$ and $q_n(x) = L_n^\mu(x)$ with $\nu \neq \mu$ but both are greater than $-1$. The second example is where $p_n(x) = L_n^\nu(z^2)$ and $q_n(x)$ is the normalized continuous dual Hahn polynomial $S_n^\mu(z^2; \sigma, \tau)$ with pure continuous spectrum (i.e., $\{\mu, \sigma, \tau\}$ are positive).

## 4. Polynomials with mixed spectra

The orthogonal polynomials $\{p_n(x)\}$ and $\{q_n(x)\}$ considered in the previous sections both have pure continuous spectra. In this section, however, we consider an interesting case with remarkable applications in nonlinear physics. The scenario is similar to that in Section 3 except that the polynomials $\{p_n(x)\}$ have a pure continuous spectrum whereas $\{q_n(x)\}$ have a mix of continuous as well as discrete spectrum. That is, the orthogonality (13) becomes

$$\int_{x_-}^{x_+} \omega(x) q_n(x) q_m(x) dx + \sum_{i=0}^{M} \chi_i q_n(x_i) q_m(x_i) = \delta_{n,m}, \tag{22}$$

where $\chi_i$ is the discrete component of the weight function and $\{x_i\}_{i=0}^{M}$ is the set of discrete spectrum of size $M+1$. We multiply both sides of Eq. (14) by $\omega(x) q_l(x)$ and integrate, calling the resulting equation, Eq. (i). Then, we write Eq. (14) at the spectrum point $x = x_i$, multiply both sides by $\chi_i q_l(x_i)$ and sum over $i = 0, 1, .., M$ calling the result Eq. (ii). Finally, we add Eq. (i) to Eq. (ii) and use the orthogonality (22) on the right side to obtain

$$\int_{x_-}^{x_+} \omega(x) p_n(x) p_m(x) q_l(x) dx + \sum_{i=0}^{M} \chi_i p_n(x_i) p_m(x_i) q_l(x_i) = d_{n,m}^l. \tag{23}$$

Since $q_0(x) = 1$, we insert $q_0(x)$ inside the integral and $q_0(x_i)$ inside the sum. Consequently, the left side of this equation becomes $[p_n(K) p_m(K)]_{0,l}$ and we obtain

$$d_{n,m}^l = [p_n(K) p_m(K)]_{0,l}, \tag{24}$$

which is identical to Eq. (19) except that the Jacobi matrix $K$ here has a mix of continuous and discrete spectra whereas the matrix $K$ in Eq. (19) has a pure continuous spectrum.

In the Appendix, we give an example where $p_n(x) = L_n^\nu(z^2)$ and $q_n(x) = S_n^\mu(z^2; \sigma, \tau)$ but with $\mu < 0$, $\mu + \sigma > 0$, and $\mu + \tau > 0$. That is, the continuous dual Hahn polynomial has a mix of discrete and continuous spectrum and the discrete spectrum is finite and of size equals to the largest integer less than $-\mu$.

## 5. Application in nonlinear physics

In this section, we use our findings above to demonstrate a remarkable phenomenon where nonlinear coupling in a physical system with a pure continuous spectrum (no bound states) can induce a dramatic change in which the spectrum becomes a mix of continuous and discrete energies (bound states). As illustration, we consider a physical system represented by the following time independent Lagrangian density (in the atomic units $\hbar = M = 1$):

$$\mathscr{L}_g(\vec{r}) = -\frac{1}{4} \psi(\vec{r}) \vec{\nabla}^2 \psi(\vec{r}) - \frac{E}{2} \psi^2(\vec{r}) - \frac{g}{3!} \vartheta(\vec{r}) \psi^3(\vec{r}), \tag{25}$$

where $g$ is a coupling parameter, $E$ is the energy and $\vartheta(\vec{r})$ is an external field. By variation of the action with respect to the wavefunction $\psi(\vec{r})$ and assuming spherical symmetry, we obtain the following wave equation for the radial component $\phi(r)$ of $\psi(\vec{r})$

$$\left[-\frac{1}{2}\frac{d^2}{dr^2} + \frac{\ell(\ell+1)}{2r^2} - E\right]\phi(r) = \frac{g}{2}\vartheta(r)\phi^2(r), \tag{26}$$

where $\ell$ is the angular momentum quantum number. For $g = 0$, the solution is well-known and could be written in terms of the spherical Bessel function. However, to facilitate our solution for $g \neq 0$, we choose an alternative, but equivalent, representation. We write the solution as an infinite expansion in a complete basis set of square integrable functions in configuration space, $\{\chi_n(r)\}_{n=0}^{\infty}$. That is, we write $\phi_0(r) = \sum_{n=0}^{\infty} A_n^0 \chi_n(r)$. All physical information about the system is contained in the expansion coefficients $\{A_n^0\}$ while the basis are generic and physically dummy. We choose the basis elements as $\chi_n(r) = (\lambda r)^{2\mu} e^{-\lambda^2 r^2/2} L_n^\nu(\lambda^2 r^2)$ where $\mu$ and $\nu$ are real dimensionless parameters such that $\mu > 0$ and $\nu > -1$. The scale parameter $\lambda$ is positive and of inverse length dimension. Thus, we write

$$\phi_0(r) = \sum_{n=0}^{\infty} A_n^0 \chi_n(r) = s^{2\mu} e^{-s^2/2} \sum_{n=0}^{\infty} A_n^0 L_n^\nu(s^2), \tag{27}$$

where $s = \lambda r$. It should be obvious that the system represented by $\phi_0(r)$ has a pure continuous energy spectrum since the corresponding Hamiltonian is the free kinetic energy operator in 3D. Substituting the ansatz (27) in Eq. (26) with $g = 0$, and using the differential equation and recursion relation of the Laguerre polynomials, we obtain $2\mu = \nu + \frac{1}{2} = \ell + 1$ and

$$z^2 A_n^0 = \left(2n + \ell + \tfrac{3}{2}\right) A_n^0 + \sqrt{n\left(n + \ell + \tfrac{1}{2}\right)} A_{n-1}^0 + \sqrt{(n+1)\left(n + \ell + \tfrac{3}{2}\right)} A_{n+1}^0, \tag{28}$$

where $z = 2\kappa/\lambda$ and $\kappa$ is the wave number (linear momentum) $\sqrt{2E}$. Comparing this recursion relation to that of the normalized Laguerre polynomial, we can write $A_n^0 \propto (-1)^n L_n^{\ell+\frac{1}{2}}(z^2)$ where the proportionality is an arbitrary function of $z$ independent of $n$ since the recursion (28) is solvable modulo such an arbitrary function. Normalization in the energy, however, fixes this arbitrary function giving

$$A_n^0(z) = (-1)^n z^{\ell+1} e^{-z^2/2} L_n^{\ell+\frac{1}{2}}(z^2). \tag{29}$$

Consequently, the free wavefunction reads as follows

$$\phi_0(r) = (zs)^{\ell+1} e^{-(z^2+s^2)/2} \sum_{n=0}^{\infty} (-1)^n L_n^{\ell+\frac{1}{2}}(z^2) L_n^{\ell+\frac{1}{2}}(s^2), \tag{30}$$

Since $\{A_n^0\}$ contain all physical information about the system, it is reassuring that they confirm that the energy spectrum associated with $\phi_0(r)$ is purely continuous since the spectrum of the Laguerre polynomial $L_n^{\ell+\frac{1}{2}}(z^2)$ is. At this point, we can follow one of two methods to obtain the solution of the wave equation (26); one is a perturbative approach and the other is non-perturbative.

## 5.1 Non-Perturbative approach:

Since the basis is complete, we can write $\phi(r) = \sum_{n=0}^{\infty} \mathcal{A}_n \chi_n(r)$ with $\chi_n(r) = s^{2\mu} e^{-s^2/2} L_n^\nu(s^2)$ and $\mathcal{A}_n \neq \mathcal{A}_n^0$, a non-perturbative solution could be obtained if the external field is chosen as $\vartheta(r) = s^{-2\mu} e^{s^2/2}$ for some positive dimensionless parameter $\mu$ because it makes the weight functions of $\phi(r)$ and of $\vartheta(r)\phi^2(r)$ on both sides of the wave equation (26) identical. Using the differential equation of the Laguerre polynomials, the wave equation becomes

$$\sum_{n=0}^{\infty} \mathcal{A}_n \left[ \left(2\mu - \nu - \tfrac{1}{2}\right) \frac{d}{dy} + \frac{(2\mu+\ell)(2\mu-\ell-1)}{4y} + \frac{y}{4} - \frac{1}{2}\left(2n + 2\mu + \tfrac{1}{2}\right) + \frac{2E}{\lambda^2} \right] L_n^\nu(y)$$

$$= \frac{g}{2} y \sum_{n=0}^{\infty} \mathcal{B}_n L_n^\nu(y) \quad (31)$$

where $y = s^2$ and $\mathcal{B}_n = \sum_{i,j=0}^{\infty} \mathcal{A}_i \mathcal{A}_j c_{i,j}^n$ with $c_{i,j}^n = [L_n^\nu(J)]_{i,j}$ and $J$ is the tridiagonal symmetric matrix (5) with $a_n = 2n + \nu + 1$ and $b_n = -\sqrt{(n+1)(n+\nu+1)}$. Now, we impose the condition that for some choice of the arbitrary parameters $\mu$ and $\nu$ we can take $2g\mathcal{B}_n = \mathcal{A}_n$. That is, in matrix notation, $\mathcal{A}_n = 2g \mathcal{A}^T c^n \mathcal{A}$ and we will have to verify the compatibility of this condition with the final solution. Substituting this back into (31) and using the differential relation of the normalized Laguerre polynomial, $y \frac{d}{dy} L_n^\nu(y) = n L_n^\nu(y) - \sqrt{n(n+\nu)} L_{n-1}^\nu$, we obtain the following with $2\mu = \nu + \tfrac{3}{2}$

$$-\tfrac{1}{4}\left(\ell + \tfrac{1}{2}\right)^2 \mathcal{A}_n = \left[(2n+\nu+1)\left(n + \tfrac{\nu}{2} + 1 - \varepsilon\right) - n - \tfrac{1}{4}(\nu+1)^2\right] \mathcal{A}_n$$

$$-\left(n + \tfrac{\nu}{2} - \varepsilon\right)\sqrt{n(n+\nu)}\, \mathcal{A}_{n-1} - \left(n + \tfrac{\nu}{2} + 1 - \varepsilon\right)\sqrt{(n+1)(n+\nu+1)}\, \mathcal{A}_{n+1} \quad (32)$$

where $\varepsilon = 2E/\lambda^2$. Comparing the recursion coefficients in this relation to those of the normalized continuous dual Hahn polynomial shown in the Appendix as (A8), we conclude that $\mathcal{A}_n(z) \propto S_n^\zeta(z^2; \sigma, \tau)$ and that

$$\sigma = \tau = \frac{\nu+1}{2}, \qquad \zeta = \frac{1}{2} - \varepsilon, \qquad z^2 = -\frac{1}{4}\left(\ell + \tfrac{1}{2}\right)^2, \quad (33)$$

Therefore, $z$ is pure imaginary for all physical scenarios, which dictates that the parameter $\zeta$ must be negative (i.e., $\varepsilon > 1/2$). This corresponds to a system with pure discrete energy spectrum, which is positive and bounded from below. In fact, using the spectrum formula of the continuous dual Han polynomial shown in the Appendix that reads $z_k^2 = -(k+\zeta)^2$, we obtain the following bound states energy spectrum

$$E_k = \frac{\lambda^2}{4}\left(2k + \ell + \tfrac{3}{2}\right), \quad (34)$$

which is that of the isotropic oscillator of frequency $\lambda^2/4$. It is worth noting that if we look at the problem as linear then we can interpret the term $-\tfrac{g}{2} \vartheta(r)\phi(r)$ in the wave equation (26) as

–8–

a potential function $V(r)$ that reads $-\frac{g}{2}\sum_{n=0}^{\infty}\mathcal{A}_n L_n^\nu(\lambda^2 r^2)$, which is not necessarily a confining potential that could be responsible for this pure discrete spectrum. Depending on the expansion coefficients $\{\mathcal{A}_n\}$, it could even be repulsive. Therefore, the phenomenon of inducing the discrete energy spectrum in the system whose original spectrum is purely continuous is due to the nonlinear coupling.

Now, we should check the validity of our initial assumption that $\mathcal{A}_n = 2g\,\mathcal{A}^T c^n \mathcal{A}$ for some proper choice of the arbitrary parameter $\nu$. In fact, even if this condition could not be satisfied then we can always find an $\ell$-independent ($z$-independent) similarity transformation matrix $\Omega$ such that $\tilde{\mathcal{A}}_n = 2g\,\mathcal{A}^T \tilde{c}^n \mathcal{A}$, where $\tilde{\mathcal{A}} = \Omega \mathcal{A}$ and $\tilde{c}^n = \sum_m \Omega_{n,m} c^m$.

## 5.2 Perturbative approach:

If we assume weak coupling then we can rewrite the wave equation (26) to first order as follows

$$\left[-\frac{1}{2}\frac{d^2}{dr^2}+\frac{\ell(\ell+1)}{2r^2}-E\right]\phi_1(r) = \frac{g}{2}\mathcal{9}(r)\phi_0^2(r). \tag{35}$$

We can also write

$$\phi_0^2(r) = s^{2\ell+2} e^{-s^2} \sum_{n,m=0}^{\infty} A_n^0 A_m^0 L_n^{\ell+\frac{1}{2}}(s^2) L_m^{\ell+\frac{1}{2}}(s^2) =$$
$$s^{2\ell+2} e^{-s^2} \sum_{n,m=0}^{\infty} \sum_{k=|n-m|}^{n+m} A_n^0 A_m^0 c_{n,m}^k L_k^{\ell+\frac{1}{2}}(s^2) = s^{2\ell+2} e^{-s^2} \sum_{k=0}^{\infty} B_k^0 L_k^{\ell+\frac{1}{2}}(s^2) \tag{36}$$

where we have used the linearization in Eq. (4) with $c_{n,m}^k = [L_k^{\ell+\frac{1}{2}}(J)]_{n,m}$ and $J$ is the tridiagonal symmetric matrix (5) with $a_n = 2n+\ell+\frac{3}{2}$ and $b_n = -\sqrt{(n+1)(n+\ell+\frac{3}{2})}$. The expansion coefficients $\{B_k^0\}$ can be written as

$$B_k^0(z) = \sum_{n,m=0}^{\infty} A_n^0 A_m^0 c_{n,m}^k = z^{2\ell+2} e^{-z^2} \sum_{n,m=0}^{\infty} \sum_{l=|n-m|}^{n+m} (-1)^{n+m} c_{n,m}^l c_{n,m}^k L_l^{\ell+\frac{1}{2}}(z^2). \tag{37}$$

Again, since $\{\chi_n(r)\}_{n=0}^{\infty}$ is a complete basis set in the radial configuration space, then we can write the solution of Eq. (35) as follows

$$\phi_1(r) = \sum_{n=0}^{\infty} A_n^1 \chi_n(r) = s^{\ell+1} e^{-s^2/2} \sum_{n=0}^{\infty} A_n^1 L_n^{\ell+\frac{1}{2}}(s^2), \tag{38}$$

where $A_n^1 \neq A_n^0$. For a feasible solution, we choose $\mathcal{9}(r) = s^{-(\ell+1)} e^{s^2/2}$ because it makes the weight functions of $\phi_1(r)$ and of $\mathcal{9}(r)\phi_0^2(r)$ on both sides of the wave equation (35) identical. Substituting $\phi_1(r)$ and $\mathcal{9}(r)\phi_0^2(r)$ in the wave equation (35), we obtain

$$z^2 A_n^1(z,g) = \left(2n+\ell+\tfrac{3}{2}\right) A_n^1(z,g) + \sqrt{n\left(n+\ell+\tfrac{1}{2}\right)} A_{n-1}^1(z,g)$$
$$+\sqrt{(n+1)\left(n+\ell+\tfrac{3}{2}\right)} A_{n+1}^1(z,g) - \frac{g}{2} B_n^0(z) \tag{39}$$

which differs from the recursion (28) for $A_n^0$ by the added last term, making its analytic solution highly nontrivial. However, if we only need the lowest few, $\{A_n^1(g,z)\}_{n=0}^N$, then we can obtain them explicitly by using the recursion (39) and knowledge of the initial value $A_0^1(z,g)$ and of $\{B_n^0(z)\}_{n=0}^N$ from (37). Alternatively, we follow an approach in which we promote the added term $-\frac{g}{2}B_n^0(z)$ to a higher intermediate order by rewriting $B_k^0$ as follows

$$B_k^0(z) = \sum_{n,m=0}^{\infty} A_n^0 A_m^0 c_{n,m}^k \mapsto B_k^{(0,1)}(z) := \sum_{n,m=0}^{\infty} A_n^0 A_m^1 c_{n,m}^k = \sum_{m=0}^{\infty} \Lambda_{k,m} A_m^1, \qquad (40)$$

where the matrix $\Lambda_{k,m} = \sum_{n=0}^{\infty} A_n^0 c_{n,m}^k = \sum_{n=0}^{\infty} c_{m,n}^k A_n^0$, where we have used the symmetry $c_{n,m}^k = c_{m,n}^k$. Consequently, Eq. (39) becomes

$$z^2 A_n^1(z,g) = \left(2n + \ell + \tfrac{3}{2}\right) A_n^1(z,g) + \sqrt{n\left(n + \ell + \tfrac{1}{2}\right)} A_{n-1}^1(z,g)$$
$$+ \sqrt{(n+1)\left(n + \ell + \tfrac{3}{2}\right)} A_{n+1}^1(z,g) - \frac{g}{2} \sum_{m=0}^{\infty} \Lambda_{n,m}(z) A_m^1(z,g) \qquad (41)$$

If the matrix $\Lambda$ can be truncated into a finite $N \times N$ dimensional matrix without any significant loss of the physics for large enough matrix size, then in matrix notation, this equation becomes

$$z^2 \begin{pmatrix} A_0^1 \\ A_1^1 \\ \times \\ \times \\ \times \\ A_{N-2}^1 \\ A_{N-1}^1 \\ A_N^1 \\ A_{N+1}^1 \\ \times \\ \times \\ \times \\ \times \end{pmatrix} = \begin{pmatrix} \times & \times & \times & \times & \times & \times & \times & & & & & & \\ \times & \times & \times & \times & \times & \times & \times & & & & & & \\ \times & \times & \times & \times & \times & \times & \times & & & & & & \\ \times & \times & \times & \times & \times & \times & \times & & & & & & \\ \times & \times & \times & \times & \times & \times & \times & & & & & & \\ \times & \times & \times & \times & \times & \times & \times & & & & & & \\ \times & \times & \times & \times & \times & \times & \times & \times & & & & & \\ & & & & & & \times & \times & \times & & & & \\ & & & & & & & \times & \times & \times & & & \\ & & & & & & & & \times & \times & \times & & \\ & & & & & & & & & \times & \times & \times & \\ & & & & & & & & & & \times & \times & \times \\ & & & & & & & & & & & \times & \times \end{pmatrix} \begin{pmatrix} A_0^1 \\ A_1^1 \\ \times \\ \times \\ \times \\ A_{N-2}^1 \\ A_{N-1}^1 \\ A_N^1 \\ A_{N+1}^1 \\ \times \\ \times \\ \times \\ \times \end{pmatrix} \qquad (42)$$

This setting makes it identical to that of the J-matrix method of scattering [15]. Consequently, we can use the tools of the J-matrix method to obtain all scattering information including phase shift, resonances and bound states.

Now for second order, the wave equation (26) is rewritten as follows

$$\left[ -\frac{1}{2}\frac{d^2}{dr^2} + \frac{\ell(\ell+1)}{2r^2} - E \right] \phi_2(r) = g\vartheta(r)\phi_1^2(r), \qquad (43)$$

where we can use (38) to write

–10–

$$\vartheta(r)\phi_1^2(r) = s^{\ell+1}e^{-s^2/2}\sum_{n,m=0}^{\infty} A_n^1 A_m^1 L_n^{\ell+\frac{1}{2}}(s^2)L_m^{\ell+\frac{1}{2}}(s^2) = s^{\ell+1}e^{-s^2/2}\sum_{n=0}^{\infty} B_n^1 L_n^{\ell+\frac{1}{2}}(s^2), \qquad (44)$$

with $B_k^1(z,g) = \sum_{n,m=0}^{\infty} A_n^1(z,g) A_m^1(z,g) c_{n,m}^k$. We write $\phi_2(r) = \sum_{n=0}^{\infty} A_n^2 \chi_n(r)$ and substitute it along with $\vartheta(r)\phi_1^2(r)$ in the second order wave equation (43). As a result, we obtain

$$\begin{aligned} z^2 A_n^2(z,g) = \left(2n+\ell+\tfrac{3}{2}\right) A_n^2(z,g) + \sqrt{n\left(n+\ell+\tfrac{1}{2}\right)} A_{n-1}^2(z,g) \\ + \sqrt{(n+1)\left(n+\ell+\tfrac{3}{2}\right)} A_{n+1}^2(z,g) - \frac{g}{2} B_n^1(z,g) \end{aligned} \qquad (45)$$

This process continues up to a desired perturbation order corresponding to a conveniently chosen accuracy. At each step $k$, we produce a set of wavefunction expansion coefficient that satisfy

$$\begin{aligned} z^2 A_n^k(z,g) = \left(2n+\ell+\tfrac{3}{2}\right) A_n^k(z,g) + \sqrt{n\left(n+\ell+\tfrac{1}{2}\right)} A_{n-1}^k(z,g) \\ + \sqrt{(n+1)\left(n+\ell+\tfrac{3}{2}\right)} A_{n+1}^k(z,g) - \frac{g}{2} B_n^{k-1}(z,g) \end{aligned} \qquad (46)$$

where $B_n^{k-1}(z,g)$ is obtained from the previous step in terms of $A_n^{k-1}(z,g)$ as

$$B_l^{k-1}(z,g) = \sum_{n,m=0}^{\infty} A_n^{k-1}(z,g) A_m^{k-1}(z,g) c_{n,m}^l. \qquad (47a)$$

Moreover, for using the J-matrix technique, we can also write

$$B_l^{(k-1,k)}(z,g) := \sum_{n,m=0}^{\infty} A_n^{k-1}(z,g) A_m^k(z,g) c_{n,m}^l. \qquad (47b)$$

## 6. Conclusion

Linearization of products of orthogonal polynomials as finite linear sum in a complete basis set of special orthogonal polynomials is an important procedure with wide range of applications in science and engineering. Most of the published linearization formulas are either too complicated to be implemented numerically, apply to specific types of polynomials, or involve integrals that are not easy to evaluate. In this work, we were able to obtain exact, very simple and general linearization formulas that are easy to implement in numerical calculations. They depend only on the recursion coefficients of the linearizing polynomials. We introduced a highly nontrivial case where the spectrum of the polynomials in the product is purely continuous whereas the linearizing polynomial has a mixed spectrum of continuous and discrete components. We presented a physically interesting application of this case in nonlinear physics. We leave for a future work the case where the polynomials have pure discrete spectra.

## **Appendix: Numerical Examples**

In this Appendix, we present some examples that illustrate and validate our findings. We start by listing all orthogonal polynomials involved in these examples along with their relevant properties. All these polynomials are written in their orthonormal version.

**Hermite:**

$$H_n(x) = \sqrt{\frac{2^n}{n!}} (x^n) \,_2F_0\left(\begin{array}{c}-n/2,-(n-1)/2\\ \underline{\quad}\end{array}\bigg| -1/x^2\right), \quad -\infty < x < +\infty. \tag{A1}$$

The recursion coefficients:

$$a_n = 0, \quad b_n = \sqrt{(n+1)/2}. \tag{A2}$$

The weight function:

$$\rho(x) = e^{-x^2}/\sqrt{\pi}. \tag{A3}$$

**Laguerre:**

$$L_n^\nu(x) = \sqrt{\frac{(\nu+1)_n}{n!}} \,_1F_1\left(\begin{array}{c}-n\\ \nu+1\end{array}\bigg| x\right), \tag{A4}$$

where $x \geq 0$, $\nu > -1$ and $(a)_n = a(a+1)(a+2)...(a+n-1) = \frac{\Gamma(n+a)}{\Gamma(a)}$ is the Pochhammer symbol, which is also known as the shifted factorial.

The recursion coefficients:

$$a_n = 2n + \nu + 1, \quad b_n = -\sqrt{(n+1)(n+\nu+1)}. \tag{A5}$$

The weight function:

$$\rho(x) = x^\nu e^{-x}/\Gamma(\nu+1). \tag{A6}$$

**Continuous dual Hahn:**

$$S_n^\mu(z^2;\sigma,\tau) = \sqrt{\frac{(\mu+\sigma)_n(\mu+\tau)_n}{n!(\sigma+\tau)_n}} \,_3F_2\left(\begin{array}{c}-n,\mu+iz,\mu-iz\\ \mu+\sigma,\mu+\tau\end{array}\bigg| 1\right). \tag{A7}$$

where $z \geq 0$. If the parameters $\{\mu,\sigma,\tau\}$ are positive then it has a pure continuous spectrum. On the other hand, if $\mu < 0$ and $\mu+\sigma > 0$, $\mu+\tau > 0$ then the spectrum is mixed and the discrete part is finite with a size equals to $M+1$ where $M$ is the largest integer less than $-\mu$.

The recursion coefficients:

$$a_n = (n+\mu+\sigma)(n+\mu+\tau) + n(n+\sigma+\tau-1) - \mu^2, \tag{A8a}$$

$$b_n = -\sqrt{(n+1)(n+\sigma+\tau)(n+\mu+\sigma)(n+\mu+\tau)}. \tag{A8b}$$

The continuous weight function:

$$\omega(z) = \frac{1}{2\pi} \frac{|\Gamma(\mu+iz)\Gamma(\sigma+iz)\Gamma(\tau+iz)/\Gamma(2iz)|^2}{\Gamma(\mu+\sigma)\Gamma(\mu+\tau)\Gamma(\sigma+\tau)}. \tag{A9}$$

The discrete weight function:

$$\chi_j = 2\frac{\Gamma(\sigma-\mu)\Gamma(\tau-\mu)}{\Gamma(\sigma+\tau)\Gamma(1-2\mu)}\frac{(-\mu-j)}{j!}\frac{(\mu+\sigma)_j(\mu+\tau)_j(2\mu)_j}{(\mu-\sigma+1)_j(\mu-\tau+1)_j}(-1)^j. \tag{A10}$$

The discrete spectrum is $z_j^2 = -(j+\mu)^2$.

## Examples for Section 2

***Example* 2.1**: $p_n(x) = q_n(x) = H_n(x)$.

The solid red trace in Figure 1a is the left side of Eq. (4) whereas the dotted blue trace is the right side with $c_{n,m}^k = [H_k(J)]_{n,m}$ where $J$ is the tridiagonal symmetric matrix of Eq. (5) with the recursion coefficients (A2). For a proper display, we rescaled the plot by multiplying both sides of Eq. (4) by the weight function (A3). In the figure, we took $n=3$ and $m=7$. On the other hand, the linearization formula in the literature reads as follows [see, for example, Eq. (3.10) in Ref. 14]

$$\hat{H}_n(x)\hat{H}_m(x) = \sum_{k=0}^{\min(n,m)} \binom{n}{k}\binom{m}{k} 2^k k! \hat{H}_{n+m-2k}(x), \tag{A11}$$

where $\hat{H}_n(x)$ is the conventional definition of the Hermite polynomial not the orthonormal version defined in (A1). That is, $\hat{H}_n(x) = \sqrt{n!/2^n}\, H_n(x)$. Figure 1b demonstrates agreement of our result with (A11) to within machine accuracy.

***Example* 2.2**: $p_n(x) = q_n(x) = S_n^\mu(z^2;\sigma,\tau)$ with the parameters $\{\mu,\sigma,\tau\}$ being positive, which makes the associate spectrum pure continuous.

The solid red trace in Figure 2 is the left side of Eq. (4) whereas the dotted blue trace is the right side with $c_{n,m}^k = [S_k^\mu(J;\sigma,\tau)]_{n,m}$ where $J$ is the tridiagonal symmetric matrix of Eq. (5) with the recursion coefficients (A8). For a proper display, we rescaled the plot by multiplying both sides of Eq. (4) by the weight function (A9). In the figure, we took $n=3$, $m=6$ and $\{\mu,\sigma,\tau\} = \{7.5, 2.0, 3.0\}$.

## Examples for Section 3

***Example* 3.1**: $p_n(x) = L_n^\nu(x)$ and $q_n(x) = L_n^\mu(x)$ with $\mu \neq \nu$.

The solid red trace in Figure 3 is the left side of Eq. (14) whereas the dotted blue trace is the right side with $d_{n,m}^k = [L_n^\nu(J)L_m^\nu(J)]_{0,k}$ where $J$ is the tridiagonal symmetric matrix of Eq. (5) with the recursion coefficients (A5) with $\nu \mapsto \mu$. For a proper display, we rescaled the plot by multiplying both sides of Eq. (14) by the weight function (A6). In the figure, we took $n=2$, $m=5$ and $\{\mu,\nu\} = \{2.0, 3.0\}$.

***Example* 3.2**: $p_n(x) = L_n^\nu(z^2)$ and $q_n(x) = S_n^\mu(z^2;\sigma,\tau)$ with the parameters $\{\mu,\sigma,\tau\}$ being positive, which makes the spectrum of $S_n^\mu(z^2;\sigma,\tau)$ pure continuous.

The solid red trace in Figure 4 is the left side of Eq. (14) whereas the dotted blue trace is the right side with $d_{n,m}^k = \left[ L_n^\nu(J) L_m^\nu(J) \right]_{0,k}$ where $J$ is the tridiagonal symmetric matrix of Eq. (5) with the recursion coefficients (A8). For a proper display, we rescaled the plot by multiplying both sides of Eq. (14) by the weight function (A6) with $x \mapsto z^2$. In the figure, we took $n = 5$, $m = 3$ and $\{\nu, \mu, \sigma, \tau\} = \{2.5, 3.5, 2.0, 1.0\}$.

**Section 4 example:** $p_n(x) = L_n^\nu(z^2)$ and $q_n(x) = S_n^\mu(z^2; \sigma, \tau)$ with the parameters $\mu < 0$, $\mu + \sigma > 0$, $\mu + \tau > 0$ making the spectrum of $S_n^\mu(z^2; \sigma, \tau)$ a mix of continuous and discrete components.

The solid red trace in Figure 5 is the left side of Eq. (14) whereas the dotted blue trace is the right side with $d_{n,m}^k = \left[ L_n^\nu(J) L_m^\nu(J) \right]_{0,k}$ where $J$ is the tridiagonal symmetric matrix of Eq. (5) with the recursion coefficients (A8). For a proper display, we rescaled the plot by multiplying both sides of Eq. (14) by the weight function (A6) with $x \mapsto z^2$. In the figure, we took $n = 2$, $m = 5$ and $\{\nu, \mu, \sigma, \tau\} = \{2.0, -5.5, 1.0 - \mu, 0.5 - \mu\}$.

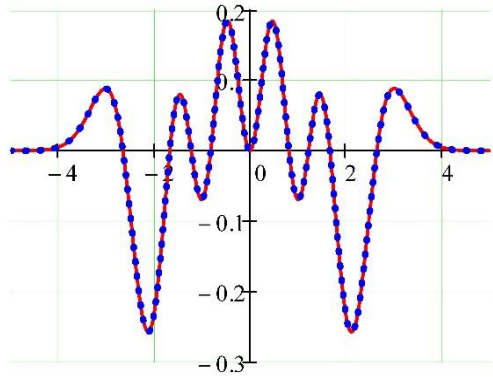

**Fig. 1a**: $p_n(x) = q_n(x) = H_n(x)$. The solid red trace is the left side of Eq. (4) whereas the dotted blue trace is the right side with $c_{n,m}^k = [H_k(J)]_{n,m}$ where $J$ is the matrix of Eq. (5) with the recursion coefficients (A2). For a proper display, we rescaled the plot by multiplying both sides of Eq. (4) by the weight function (A3). In the figure, we took $n=3$ and $m=7$.

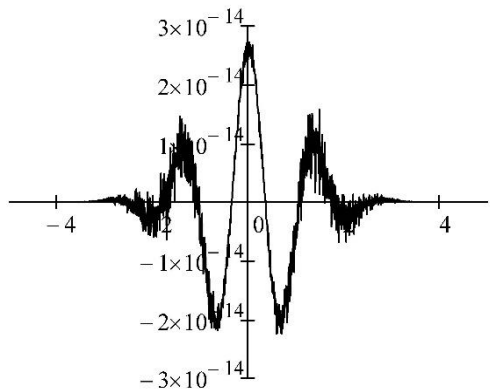

**Fig. 1b**: Difference between our result shown in Fig. 1a and that of formula (A11) to within machine accuracy.

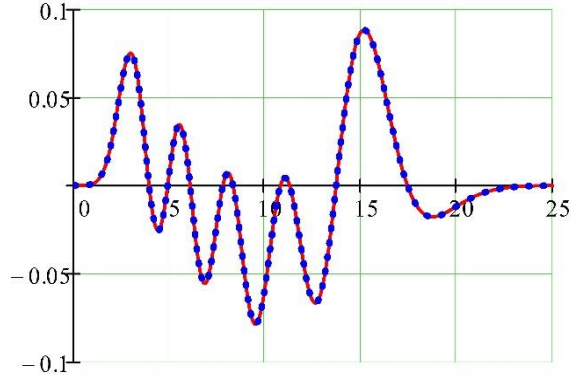

**Fig.2**: $p_n(x) = q_n(x) = S_n^\mu(z^2;\sigma,\tau)$ with $\mu > 0$. The solid red trace is the left side of Eq. (4) whereas the dotted blue trace is the right side with $c_{n,m}^k = \left[ S_k^\mu(J;\sigma,\tau) \right]_{n,m}$ where $J$ is the matrix of Eq. (5) with the recursion coefficients (A8). For a proper display, we rescaled the plot by multiplying both sides of Eq. (4) by the weight function (A9). In the figure, we took $n = 3$, $m = 6$ and $\{\mu,\sigma,\tau\} = \{7.5, 2.0, 3.0\}$.

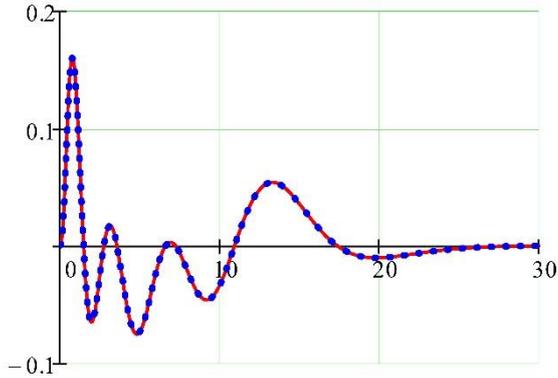

**Fig. 3**: $p_n(x) = L_n^\nu(x)$ and $q_n(x) = L_n^\mu(x)$ with $\mu \neq \nu$. The solid red trace is the left side of Eq. (14) whereas the dotted blue trace is the right side with $d_{n,m}^k = \left[ L_n^\nu(J) L_m^\nu(J) \right]_{0,k}$ where $J$ is the matrix of Eq. (5) with the recursion coefficients (A5) with $\nu \mapsto \mu$. For a proper display, we rescaled the plot by multiplying both sides of Eq. (14) by the weight function (A6). In the figure, we took $n = 2$, $m = 5$ and $\{\mu,\nu\} = \{2.0, 3.0\}$

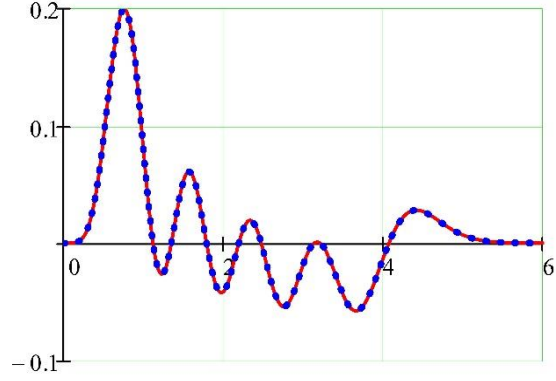

**Fig. 4**: $p_n(x) = L_n^\nu(z^2)$ and $q_n(x) = S_n^\mu(z^2;\sigma,\tau)$ with $\mu > 0$. The solid red trace is the left side of Eq. (14) whereas the dotted blue trace is the right side with $d_{n,m}^k = \left[ L_n^\nu(J) L_m^\nu(J) \right]_{0,k}$ where $J$ is the matrix of Eq. (5) with the recursion coefficients (A8). For a proper display, we rescaled the plot by multiplying both sides of Eq. (14) by the weight function (A6) with $x \mapsto z^2$. In the figure, we took $n = 5$, $m = 3$ and $\{\nu,\mu,\sigma,\tau\} = \{2.5, 3.5, 2.0, 1.0\}$.

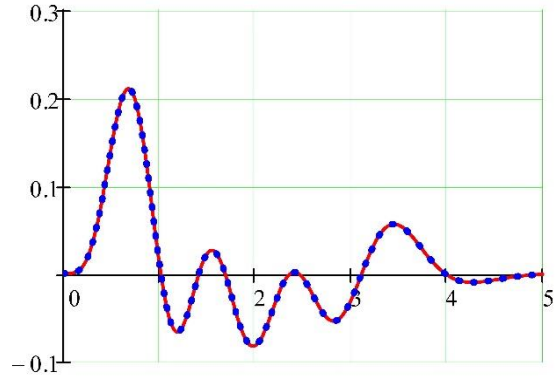

**Fig. 5**: $p_n(x) = L_n^\nu(z^2)$ and $q_n(x) = S_n^\mu(z^2;\sigma,\tau)$ with $\mu < 0$. The solid red trace is the left side of Eq. (14) whereas the dotted blue trace is the right side with $d_{n,m}^k = \left[ L_n^\nu(J) L_m^\nu(J) \right]_{0,k}$ where $J$ is the matrix of Eq. (5) with the recursion coefficients (A8). For a proper display, we rescaled the plot by multiplying both sides of Eq. (14) by the weight function (A6) with $x \mapsto z^2$. In the figure, we took $n = 2$, $m = 5$ and $\{\nu,\mu,\sigma,\tau\} = \{2.0, -5.5, 1.0-\mu, 0.5-\mu\}$.